\documentclass[oneside,11pt]{article}
\usepackage{epsfig, amsfonts, amsmath, amsthm,amsbsy}
\usepackage{color}
\usepackage{subfigure}

\newtheorem{theorem}{Theorem}
\newtheorem{lemma}{Lemma}

\newtheorem{remark}{Remark}

\newtheorem*{example-non}{Example}

\def\ve{\varepsilon}

\begin{document}

\title{Fitted finite element methods for singularly perturbed elliptic  problems of convection-diffusion type}
\author{A. F. Hegarty 
\thanks{Department of Mathematics and Statistics, University of Limerick, Ireland.} \and E.\ O'Riordan \thanks{ School of Mathematical Sciences, Dublin City University, Dublin 9, Ireland.} }

\maketitle

\begin{abstract}

 Fitted finite element methods are constructed for a singularly perturbed convection-diffusion problem in two space dimensions.  
Exponential splines as basis functions are combined with Shishkin meshes to obtain a stable parameter-uniform numerical method.
These schemes satisfy a discrete maximum principle. 
In the classical case, the numerical approximations converge,  in $L_\infty$, at a rate of second order and the approximations converge at a rate of first order for all values of the singular perturbation parameter.
\bigskip


\end{abstract}


\section{Introduction}

The numerical solution of singularly perturbed convection-diffusion problems  of the form
\[
-\ve \triangle u + \vec a \cdot \nabla u + bu =f, \quad 0 < \ve \leq 1,
\]
 presents several computational difficulties. A  key issue is to avoid 
spurious oscillations in the numerical approximations. These oscillations can be damped using various variants of streamline diffusion finite element methods (SDFEM) \cite{john1, john2}.
However, within a finite element framework it is difficult to generate a numerical method that preserves the inverse-monotonicity of the differential operator \cite{gabriel} and guarantees no spurious oscillations in the numerical approximations. A significant reduction of these oscillations can be achieved
by astute choices of stabilization parameters, but small oscillations can still be problematic  when solving a system of nonlinear  partial differential equations if  one   equation within the system  is a singularly perturbed differential equation of convection-diffusion type. 
In this paper, we choose to only consider discretizations which  are inverse monotone. Moreover, we are solely interested in parameter-uniform \cite{fhmos} numerical methods. That is, numerical methods for which an error bound on the numerical approximations $U$ of the form
\[
\Vert u - U \Vert_\infty \leq C N^{-p}, \quad p >0,
\]
can be established. Here $N$ is the number of elements used in any coordinate direction, $\Vert \cdot \Vert _\infty$ is the pointwise $L_\infty$ norm and the  error constant $C$ (used throughout this paper) is independent of $N$ and $\ve$. 

Within a  finite difference (or finite volume) framework, inverse monotonicity can be retained by employing standard upwinding or variants of upwinding. However,  upwinding  limits the order of convergence of the numerical method to first order. 
For non-singularly perturbed problems with smooth solutions, second order convergence is easily achieved using central difference schemes or classical Galerkin with bilinear basis functions. For singularly perturbed problems, it is desirable that any parameter-uniform numerical method designed to 
be inverse-monotone would also be second order  when the inverse of the singular perturbation parameter was not  large compared  to the dimensions of the discrete problem (i.e.,  if $\ve^{-1} \leq CN$).
In the case of singularly perturbed ordinary differential equations, fitted operator methods (which are nodally exact schemes in the case of  constant coefficients)   have this property at the nodes of a  uniform grid. These fitted operator methods can be generated within a finite element framework by incorporating a tensor product of one dimensional 
exponential $L$-splines or $L^*$-splines into the trial or test space \cite{rst2}. Hemker \cite{hemker} was the first to examine these exponential basis functions.  In this paper we will use these exponential splines in our choice of  trial and test space.

If one employs a uniform mesh with a tensor product of one dimensional  $L$-splines as basis functions in a Galerkin finite element method,  then one has \cite{orst4}, 
\[
\Vert u-U \Vert _E \leq  C N^{-1/2},
\]
where $\Vert \cdot \Vert _E$ is the standard $\ve$-weighted energy norm.  For the same choice of basis functions on a uniform rectangular mesh,  D\"{o}rfler \cite{dorfler} established  error bounds in a range of $L_p$-norms. The numerical performance in the $L_2$ and $L_\infty$ norms, of different combinations of exponential basis functions within a Petrov-Galerkin framework (where the trial and test space are not necessarily the same), on a uniform mesh,  was examined in \cite{hegorst}. However, as noted in \cite[Remark 3.7]{hegorst}, these schemes are not stable.

In \cite{styor4}, it is established that using  bilinear basis functions on a Shishkin mesh \cite{fhmos,mos2}  in a Galerkin framework, yields
\[ \Vert u-U \Vert _E \leq  C N^{-1} \ln N \]
and in \cite{zhimin} the superclose result \[ \Vert u_I-U \Vert _E \leq C \ve N^{-3/2} + C (N^{-1} \ln N)^2 \] 
is established, where $u_I$ is the interpolant of the continuous solution  $u$ in the trial space. 
In \cite{Martin+Lutz2003}, the same superclose bound is established  for  bilinear SDFEM on a Shishkin mesh.
 For piecewise bilinear  SDFEM on a Shishkin mesh \cite[pp 391]{rst2}, one has
\[ \Vert u_I-U \Vert _{SD} \leq C \ve N^{-3/2}+ C N^{-2} (\ln N)^{2}, \]
where $\Vert \cdot \Vert _{SD}$ is the steamline diffusion norm \cite[pp 303]{rst2}.
Using an improved pointwise interpolation bounds \cite{zhimin} and this supercloseness result, with a finite difference argument in the corner area, one can get the following result \cite[pp 400]{rst2}, for piecewise bilinear  SDFEM on a Shishkin mesh,
\[ \Vert u_I-U \Vert _\infty \leq C \ve N^{-1/2}+ C N^{-1} (\ln N)^{2}; \]
with higher orders established within the fine mesh regions. However, none of these schemes are guaranteed to be inverse-monotone for all values of the singular perturbation parameter. In this paper, given that we design the numerical methods to be inverse-monotone, we establish estimates in the pointwise $L_\infty$ norm. 
This norm  identifies all boundary  and corner layer  functions that can exist in the solution of singularly perturbed problems. In addition, we employ a layer adapted mesh (of Shishkin type) which means that we will achieve parameter-uniform bounds on the global pointwise error.

In this paper, we combine the benefits of a Shishkin mesh  with exponential splines as basis functions. Using a Shishkin mesh with basic upwinding, 
yields a globally pointwise accurate numerical approximation, satisfying a parameter-uniform bound of the form \cite{hemkera}
\[
\Vert u- U \Vert _\infty \leq C N^{-1} \log N.
\]
The log defect in this error bound can be removed if one uses upwinding and a Bakhvalov mesh \cite{nhan}, instead of a Shishkin mesh. 
Below in the later sections,  we use a tensor product of exponential  basis functions on a Shishkin mesh in a finite element formulation.  Within this framework, we establish the first order global $L_\infty$-norm error  bound
\[
\Vert u- U \Vert _\infty \leq C N^{-1}
\]
and  second order  in the classical case of $\ve=1$. Moreover,  the numerical schemes satisfy a discrete maximum principle. 

\section{Continuous problem}
Consider the singularly perturbed convection-diffusion elliptic problem 
\begin{subequations}\label{prob-2d-elliptic}
\begin{eqnarray}
 Lu \equiv -\ve \triangle u+ \vec a \cdot \nabla u+bu=f, \quad (x,y) \in \Omega :=(0,1)^2 ;\\
u=0, \quad (x,y) \in \partial \Omega ;\\
\vec a =(a_1,a_2), \quad a_1 > \alpha _1 >0,\ a_2> \alpha _2 > 0; \quad b \geq 0; \\
f \in C^{1,\gamma }(\overline{ \Omega }), \ f(0,0)=f(1,0)=f(0,1)=f(1,1)=0 \label{basic-comp}.
\end{eqnarray}

 The remaining data $a_1,a_2,b$ are assumed to be sufficiently regular so that $u \in C^{3,\gamma
}(\overline{ \Omega })$ and such that only
exponential boundary layers appear near the outflow edges $x=1,\
y=1$ and a simple corner layer appears in the vicinity of $(1,1)$.
This corner layer is induced not by any lack of sufficient
compatibility, but by the presence of the singular perturbation
parameter. In this case, there is no loss in generality in dealing
with homogeneous boundary data. 

The
solution $u$  can be  decomposed \cite{mos2}   into the sum
\[ u=v+w_T+w_R+w_{TR}.
\] Here $v$ is the regular component, $w_R$ is a regular boundary
layer function associated with the  edge $x=1$, $w_T$ is a
regular boundary layer function associated with the edge
 $y=1$ and $w_{TR}$ is a corner layer function associated with the
 corner $(1,1)$. The decomposition into regular and layer
components is defined \cite{mos2, hemkera}, so that each of the layer functions
satisfy the homogeneous differential equation $Lw=0$.

By assuming the additional regularity and compatibility conditions  \cite{hemkera}
\begin{eqnarray}\label{extra-reg}
a_1,a_2,b,,f \in C^{5,\gamma }(\overline{ \Omega }), \quad \frac{\partial ^{i+j} f}{\partial x^i \partial y ^j} (0,0)=0,\ 0 \leq i+j \leq 4
\end{eqnarray}
\end{subequations}
 on the data,  one can establish  the following bounds on the regular component 
\begin{equation}\label{regular-bnd}
\Bigl \Vert \frac{\partial ^{i+j}v}{ \partial  x^i\partial  y^j}  \Bigr\Vert _\infty  \leq
C (1+\ve ^{2-i-j}), \ i+j \leq 3.
\end{equation}

 The regular layer component $w_R$ is the solution of   the problem
  \begin{subequations}
\begin{eqnarray}
Lw_R=0, \quad (x,y) \in \Omega , \\
w_R(1,y) = (u-v) (1,y), \quad w_R(0,y)=w_L(x,0)=0, \\
w_R(x,1)=g(x).
\end{eqnarray}\end{subequations} Using a maximum principle we can deduce that
\begin{subequations}\label{edge-bnd2}
\begin{equation}\label{edge-bnd}
\vert w_R (x,y) \vert \leq C e^{-\frac{\alpha _1(1-x)}{\ve}},
\quad (x,y) \in \Omega .
\end{equation}
Using the arguments in   \cite[Chapter 12]{mos2}, coupled with the local bounds given in \cite[132--134]{lady}  and the arguments in \cite{ria}, one can deduce that 
\begin{equation}\label{orth-bnd}
\Bigl \vert \frac{\partial ^{i+j}w_R}{ \partial  x^i\partial  y^j} (x,y) \Bigr\vert  \leq
C\ve ^{-i}(1+\ve ^{1-j})e^{-\frac{\alpha _1(1-x)}{\ve}}, \ i+j \leq 3, \ (x,y) \in \Omega.
\end{equation}
Corresponding bounds hold for $w_T$. 
 Finally, we consider the corner layer function, which is defined  as
follows:
\begin{eqnarray}
Lw_{TR}=0,  \quad (x,y) \in \Omega , \\
w_{TR}(x,0) =   w_{TR}(0,y)= 0,\\
w_{TR}(1,y)=-w_T(1,y), \   w_{TR}(x,1)=-w_R(x,1);
\end{eqnarray} and we have the bounds
\begin{equation}\label{corner-bnd}
\vert w_{TR}(x,y) \vert \leq Ce^{-\frac{\alpha _1(1-x)
}{\ve}}e^{-\frac{\alpha _2 (1-y)}{\ve}}, \quad (x,y) \in  \Omega .
\end{equation}
Repeating the argument that led to (\ref{orth-bnd}), we deduce the bounds
\begin{equation}\label{corner-bnds}
\Bigl \vert \frac{\partial ^{i+j}w_{TR}}{ \partial  x^i\partial  y^j} (x,y) \Bigr \vert  \leq C\ve ^{-(i+j)}e^{-\frac{\alpha _1(1-x)
}{\ve}}e^{-\frac{\alpha _2 (1-y)}{\ve}}, \quad
i+j \leq 3, \ (x,y) \in \Omega .\end{equation}\end{subequations}

\begin{remark}
Andreev \cite{andreev} establishes this decomposition of the solution and derives  bounds on the regular and layer components, while only imposing the  compatibility constraints (\ref{basic-comp}), (\ref{extra-reg}) at the inflow corner $(0,0)$. No compatibility is imposed at the other three corners. 
However, in order to avoid dealing with the additional corner singularities at the corners $(0,1),(1,0), (1,1)$, we confine our analysis to the case where, in addition to (\ref{extra-reg}) at the inflow corner, the basic compatibility  (\ref{basic-comp}) is assumed to hold at all four corners. 
\end{remark}

\section{Finite element framework}
A weak form of problem (\ref{prob-2d-elliptic}) is: find $u\in H_0^1(\Omega)$ such that
\begin{subequations}\label{weak-form}
\begin{eqnarray}
B(u,v) = (f,v), \ v \in H^1_0(\Omega),\qquad \hbox{where}\\ 
B(u,v):= (\ve \nabla u, \nabla v)+(\vec a\cdot \nabla u, v)+(bu, v);
\end{eqnarray}
\end{subequations}
where $(u,v)$ is the standard inner product in $L^2$ .

The domain is discretized $\bar \Omega = \cup _ {i,j=1}^N \bar \Omega _{i,j}$ by the rectangular elements
\[
\bar \Omega _{i,j} := [x_{i-1}, x_i] \times [y_{j-1} y_j ], 1 \leq i,j \leq N;
\]
where the nodal points are given by the following sets
\[
\bar \omega _x := \{ x_i \vert x_i=x_{i-1} +h_i \} _{i=0}^N, \quad \bar\omega _y := \{ y_j \vert y_j=y_{j-1} +k_j \} _{j=0}^N.
\]
We define the average mesh steps with
\[
\bar h_i := \frac{h_{i+1}+h_i}{2}, \qquad \bar k_j:= \frac{k_{j+1}+k_j}{2}.
\]
This mesh  is a tensor
product  of two piecewise-uniform one dimensional Shishkin meshes \cite{mos2}. The mesh $\bar \omega _x$ places $N/2$ elements into both $[0,1-\tau _x]$ and $[1-\tau _x,1]$, where the transition parameters are taken to be
 \begin{equation}\label{transition parameters}
\tau _x = \min \{ 0.5, 2\frac{\ve}{\alpha _1} \ln N \} \quad \hbox{
and (analogously for  $\bar \omega _y$)} \quad \tau _y = \min \{ 0.5, 2\frac{\ve}{\alpha _2} \ln N
\} .
 \end{equation}
We denote the set of nodal points in this Shishkin mesh by
\[
\Omega ^N := \omega _x \times \omega _y.
\]

The trial and test space  will be denoted by $S^N,T^N \subset H_0^1(\Omega)$, respectively.
The choice of trial and test functions $\{ \phi_{i,j} (x,y) :=\phi _i(x)\phi^j(y) \} _{i=1}^{N-1} \in S^N, \{ \psi_{i,j} (x,y) :=\psi _i(x)\psi^j(y) \} _{i=1}^{N-1}\in T^N$ are simply a tensor product of  one dimensional functions, with the following standard properties
\begin{eqnarray*}
0 \leq \phi _i(x) \leq 1, \quad \phi _i (x_j) = \delta _{i,j}, \\ \hbox{supp}(\phi _i(x) ) = (x_{i-1}, x_{i+1}), \quad \phi _i(x) + \phi _{i-1}(x) =1,\ x \in  (x_{i-1}, x_{i}).
\end{eqnarray*}
Observe that we denote trial and test functions in the horizontal direction with subscripts and in the vertical direction with superscripts. 

An approximate solution $U\in S^N$ to the solution of  problem (\ref{weak-form}) is:
 find $U\in S^N$ such that
\begin{equation}
B(U,V) = (f,V), \quad \forall  V \in T^N.
\end{equation}
We denote the nodal values $U(x_i,y_j)$ simply by $U_{i,j}$. Hence
\[
U (x,y) = \sum _{i,j=1}^N U_{i,j} \phi _i(x)\phi^j(y).
\]
We now define $\bar a_{1,i}$ as piecewise constant functions, which approximate the convective coefficient $a_1$ by constant values on each element $\bar \Omega _{i,j}$. 
For example, one possible choice would be
\begin{eqnarray*}
\bar a_1 (x,y) = \bar a_{1,i} := \frac{a_1(x_{i-1},y_j)+a_1(x_{i},y_j)}{2}, \  (x,y)\in (x_{i-1},x_i] \times (y_{j-1},y_j],\\  \bar a_1 (0,y) =\bar a_1 (x_1,y), \ y \in [0,1]; \quad  \bar a_1 (x,0) =\bar a_1(x,y_1), \ x \in (0,1];\\ \\
\bar a_2 (x,y) = \bar a_{2,j} := \frac{a_2(x_i,y_{j-1})+a_2(x_i,y_j)}{2}, \ (x,y)\in   (x_{i-1},x_i]  \times (y_{j-1},y_j], \\ \bar a_2 (0,y) =\bar a_2 (x_1,y), \ y \in (0,1]; \quad  \bar a_2 (x,0) =\bar a_2(x,y_1), , \ x \in [0,1]; \\ \\
\bar f(x,y) = \bar f_{i,j}, \quad  (x,y)\in (x_{i-1},x_i] \times (y_{j-1}, y_j], \
\bar f(x,y) = \bar f_{1,1}, \ (x,y) \in \bar \Omega _{1,1} \\\hbox{where} \quad 
\bar f_{i,j} := \frac{f(x_{i-1},y_{j-1})+f(x_{i-1},y_{j})+f(x_{i},y_{j-1})+f(x_{i},y_{j})}{4};\\ \\
\bar b(x,y) = b(x_i,y_j) , \ (x,y)\in (x_{i-1},x_i] \times (y_{j-1}, y_j], \\
\bar b(x,y) = b(x_1,y_1), \ (x,y) \in \bar \Omega _{1,1}. 
\end{eqnarray*}
Approximating the data $a_1, a_2,b, f$ in the weak form by the piecewise constant functions $\bar a, \bar b,  \bar f$ means that  all integrals (in the weak form) can  be evaluated exactly.

In addition, we will lump all zero order terms, which yields increased stability and gives a simpler structure to the definition of the system matrix. That is, we introduce the additional quadrature rules
\begin{eqnarray*}
(a_1\frac{\partial \phi _{n,m}}{\partial x} , \psi _{i,j}) &\approx& (\bar a_1\frac{\partial \phi _{n}}{\partial x}(x), \psi _i(x)) (1, \psi^{j}(y)) \delta _{m,j}, \\
(a_2\frac{\partial \phi _{n,m}}{\partial y} , \psi _{i,j})&\approx& (\bar a_2\frac{\partial \phi ^{m}}{\partial y}, \psi ^j(y)) (1, \psi_{i}(x)) \delta _{n,i},\\
(b\phi _{n,m}, \psi _{i,j} ) &\approx& b(x_i,y_j) (1, \psi_{i}(x)) (1, \psi^{j}(y))  \delta _{n,i} \delta _{m,j},
\end{eqnarray*}
where $\delta _{i,j}$ is the Kronecker delta. Then we have the following quadrature rule
\begin{eqnarray*}
(\bar b U,V) := \sum _{i,j=1}^N  b(x_i,y_j) U(x_i,y_j)V(x_i,y_j)(1, \psi_{i}(x)) (1, \psi^{j}(y)); \\
(\bar a_1 U_x,V) := \sum _{j=1}^N(\bar a_1 U_x(x,y_j), V(x,y_j)) (1, \psi^{j}(y)); \\
(\bar a_2 U_y,V) := \sum _{i=1}^N   (\bar a_2 U_y(x_i,y),V(x_i,y)) (1, \psi_{i}(x)).
\end{eqnarray*}
The associated discrete weak problem is: find $U := \sum _{i,j=1}^N U_{i,j} \phi _i(x) \phi ^j(y) \in S^N$ such that  
\begin{eqnarray}\label{discrete-weak-form}
\bar B(U,\psi_{i,j})= (\bar f,\psi_{i,j}), \quad \forall \psi_{i,j}\in T^N; \\
\bar B(U,V):= \ve (U_x, V_x)  +\ve (U_y, V_y) +  (\bar a_1U_x, V)  + (\bar a_2 U_y,V) + (\bar b U,V). \nonumber
\end{eqnarray}

The associated finite difference scheme to this finite element method is:
\begin{subequations}\label{Stencil}
\begin{equation}
\sum _{n=i-1}^{i+1} \sum _{m=j-1}^{j+1} \alpha _{n,m}  U_{n,m} + b(x_i,y_j) (1, \psi_{i,j}) U_{i,j} =  \sum _{n=i-1}^{i+1} \sum _{m=j-1}^{j+1} \gamma _{n,m} f(x_n,y_m)
\end{equation}
where the coefficients are given by:
\begin{eqnarray} 
\alpha _{i-1,j+1}=\alpha _{i+1,j+1}=\alpha _{i1,j-1}=\alpha _{i-1,j-1}=0 \label{10b}\\
\alpha _{i-1,j}=R^-_x(1,\psi ^j(y))  \quad \alpha _{i+1,j}=R^+_x(1,\psi ^j(y); \\
\alpha _{i,j-1}=R^-_y(1,\psi _i(x))\quad \alpha _{i,j+1}=R^+_y(1,\psi _i(x))   \\
\alpha _{i,i}=-(R^-_x+R_x^+) (1,\psi ^j(y)) -(R^-_y+R_x^y)(1,\psi _i(x)) 
  \\  \nonumber
\left(\begin{array}{ccc}\gamma _{i-1,j+1}&\gamma _{i,j+1}&\gamma _{i+1,j+1} \\ \gamma _{i-1,j}&\gamma _{i,j}&\gamma _{i+1,j} \\\gamma _{i-1,j-1}&\gamma _{i,j-1}&\gamma _{i+1,j-1}\end{array}\right) =\frac{1}{4}
 \left(\begin{array}{c}Q_y^+ \\ Q_y^C \\ Q_y^- \end{array}\right). \left(\begin{array}{ccc}Q_x^- \ Q_x^C \ Q_x^+ \end{array}\right) .
\end{eqnarray}
The off-diagonal matrix entries in (\ref{10b}) are zero due to our use of lumping. The elements in $R^\pm_x,$ and $Q^\pm_x$ are defined by
\begin{eqnarray*}
R_x^-:= (\phi '_{i-1}, \ve \psi '_i +\bar a_{1} \psi _i ), \  R_x^+:= (\phi '_{i+1}, \ve \psi '_i +\bar a_{1} \psi _i), \\ R_x^c:=(\phi '_{i}, \ve \psi '_i +\bar a_{1} \psi _i) =-R^-_x-R^+_x, \quad  Q_x^C:= \int _{x=x_{i-1}}^{x_{i+1}}\psi _i (x) \ dx, = Q_x^-+Q_x^+, \\ 
Q^-_x :=  \int _{x=x_{i-1}}^{x_i} \psi _i (x) \ dx, \ Q^+_x := \int _{x=x_{i}}^{x_{i+1}} \psi _i (x) \ dx .
\end{eqnarray*}
The elements in $R^\pm_y$ and $Q^\pm_y$ are defined analogously. 
\end{subequations}

We now introduce the unit $L$-spline $B(t;\rho)$ and the unit $L^*$-spline $B^*(t;\rho)$, where $\rho$ is a positive constant,  as the solutions of the two point boundary value problem
\begin{eqnarray*}
-B_{tt} +\rho B_t =0,\ t \in (0,1); \quad B(0)=0,\ B(1) =1;\\
B^*_{tt} +\rho B^*_t =0,\ t \in (0,1); \quad B^*(0)=0,\ B^*(1) =1.
\end{eqnarray*}
That is
\[
B(t) = \frac{e^{-\rho(1-t) }-e^{-\rho }}{1-e^{-\rho}} \quad \hbox{and} \quad B^*(t) = \frac{1-e^{-\rho t}}{1-e^{-\rho}}.
\]
Observe that the derivatives at the end-points of the interval $(0,1)$ are given by
\[
B_t(1)=B^*_t(0) = \sigma (\rho ) ;\  B^*_t(1)=B_t(0) = \sigma (-\rho); \quad \sigma (x) := \frac{x}{1-e^{-x}}. 
\]
Associated with these unit splines, we define the set of $\bar L$-spline basis functions  (denoted by $B_{i,j}(x,y)$) and the the set of $\bar L^*$-spline basis functions (denoted by  $B^*_{i,j}(x,y)$) as follows:
over each computational cell $[x_{i-1},x_{i+1}] \times [y_{j-1},y_{j+1}]$ 
\begin{eqnarray*}
B_{i,j}(x,y) := \Phi _{i,j} (x) \Phi ^{i,j} (y); \quad \rho _{1:i,j} := \frac{\bar a_1(x_i,y_j) h_i}{\ve} ;\ \rho _{2:i,j} := \frac{\bar a_2(x_i,y_j) k_j}{\ve}  \\
\Phi_{i,j}(x) := \Bigl\{ \begin{array}{ll} B(\frac{x-x_{i-1}}{h_i}; \rho _{1:i,j}), \ x_{i-1} \leq x \leq x_i \\ 
1- B(\frac{x-x_{i}}{h_{i+1}}; \rho _{1:i+1,j}), \ x_i \leq x \leq x_{i+1} \end{array}  \\
\Phi^{i,j}(x) := \Bigl\{ \begin{array}{ll} B(\frac{y-y_{j-1}}{k_j}; \rho _{2:i,j}), \ y_{j-1} \leq y \leq y_j\\ 
1- B(\frac{y-y_j}{k_{j+1}}; \rho _{2:i,j+1}), \ y_j \leq y \leq y_{j+1} \end{array}  \\
\end{eqnarray*}
 The basis functions $B^*_{i,j}(x,y)$ are defined analogously. 

If we use either $\bar L$-splines in the trial space and any choice of test space or  $\bar L^*$-splines in the test space and any choice of trial space,
then in either case  we have that
\begin{subequations}\label{def-R}
\begin{eqnarray}
R_x^- = -\frac{\ve \sigma (\rho _{1;i,j})} {h_i}, \  R_x^+= -\frac{\ve \sigma (-\rho _{1;i+1,j})}{h_{i+1}}, \\
R_y^- = -\frac{\ve \sigma (\rho _{2;i,j})}{k_j}, \  R_y^+= -\frac{\ve \sigma (-\rho _{2;i+1,j})}{k_{j+1}}, 
 \end{eqnarray}
\end{subequations}

Combined with (\ref{def-R}), we consider the following three possible choices for the remaining terms in (\ref{Stencil}). 

 \begin{subequations}\label{Fitted-variants-lumped}
\begin{enumerate}
 \item $\bar L^*$ splines in the test space and bilinear trial functions. 
\begin{eqnarray}\label{test_L^*_lumped}
Q^-_x= h_i\frac{\sigma (\rho _{1;i,j}) -1}{\rho _{1;i,j} },\quad  Q^+_x=  h_{i+1}\frac{1-\sigma ( -\rho _{1;i+1,j})}{\rho _{1;i+1,j}},
 \end{eqnarray}
\item $\bar L$ splines in the test space and bilinear trial functions. 
\begin{eqnarray}\label{test_L_lumped}
 Q^-_x= h_i\frac{1-\sigma (-\rho _{1;i,j}))}{\rho _{1;i,j}},\quad   Q^+_x=  h_{i+1}\frac{\sigma (\rho _{1;i+1,j})-1}{\rho _{1;i+1,j}};\  
 \end{eqnarray}
 \item  Bilinear test functions and $\bar L$ splines in the trial space.
\begin{eqnarray}\label{test_hat_lumped}
 Q^-_x= \frac{h_i}{2},\ Q^+_x= \frac{ h_{i+1}}{2}
 \end{eqnarray}
\end{enumerate}
\end{subequations}
with analogous definitions for $Q^+_y,Q^-_y$ in each case. All of these fitted schemes have an M-matrix structure and we, hence, have guaranteed stability.

\begin{remark}
We can write these  fitted schemes in finite difference notation as follows
\begin{eqnarray*}
\frac{Q^C_y}{\bar h_i \bar k_j} \bigl (-\ve h_{i+1}D^+_x(\sigma (-\rho _{1:i,j}) D^-_x)+\bar a_1 h_{i}D^-_x\bigr) U(x_i,y_j) \\
+ \frac{Q^C_x}{\bar h_i \bar k_j} \bigl (-\ve k_{j+1}D^+_y(\sigma (-\rho _{2:i,j}) D^-_y+\bar a_2 k_jD^-_y\bigr) U(x_i,y_j)\\
+ b(x_i,y_j) \frac{Q^C_xQ^C_y}{\bar h_i \bar k_j}  U(x_i,y_j)= \frac{1}{\bar h_i \bar k_j}   \sum _{n=i-1}^{i+1} \sum _{m=j-1}^{j+1} \gamma _{n,m} f(x_n,y_m) .
 \end{eqnarray*}
In the case of constant coefficients $a_1(x,y)= \alpha _1, a_2(x,y)=\alpha _2$ and $b(x,y) \equiv 0$, this finite difference scheme is exact for the boundary layer functions
$e^{-\alpha _1(1-x)/ \ve}, e^{-\alpha _2(1-y)/ \ve} $ and the corner layer function $e^{-\alpha _1(1-x)/ \ve}e^{-\alpha _2(1-y)/ \ve}$ on an arbitrary mesh.
\end{remark}
\begin{remark}\label{remark2}
Lin\ss \  \cite{linss2}  examined a class of fitted finite difference operators (arising from a finite volume formulation) on a tensor product of Shishkin meshes. Using the $(L_\infty, L_1)$ stability argument developed by Andreev \cite{andreev2003}, Lin\ss \ established that for $u \in C^4(\bar \Omega) $, then
\[
\vert (u-U)(x_i,y_j) \vert  \leq  C N^{-1},
\]
if the fitting factor satisfies certain conditions (see \cite[(3), pg. 248]{linss2}).
If we formally set $Q_x^-=Q_x^+=Q_y^+=Q_y^+ =0, Q_x^C= \bar h_i, Q_y^C= \bar k_j$ in  (\ref{Fitted-variants-lumped}), the resulting finite difference scheme scheme fits into the framework of fitted finite difference schemes analysed in Lin\ss \  \cite{linss2}.
\end{remark}

\begin{remark} The above numerical schemes can also be applied to the problem
\[
-\ve \triangle u + a_1u_x + bu =f.
\]
The solution will now have a regular layer near the outflow boundary $x=1$, two characteristic layers along the sides $y=0$ and $y=1$ and, assuming sufficient compatibility at the inflow corners, corner layers in neighbourhoods of the two outflow corners $(1,0), (1,1)$. An appropriate 
piecewise uniform Shishkin mesh can be constructed to capture these layers \cite{hemkerb}. For this problem, some of the terms in the above fitted schemes simplify to
\[
R^-_y = - \frac{\ve}{k_j}, R^+_y = - \frac{\ve}{k_{j+1}}, \quad Q^-_y =  \frac{k_j}{2}, Q^+_y = \frac{k_{j+1}}{2}.
\]
\end{remark}

\section{Error analysis in $L_\infty$}

To identify the truncation errors associated with each coordinate direction, we introduce the following notation for one dimensional differential operators
\begin{eqnarray*}
Lu = L_xu+L_yu + bu, \\ 
\hbox{where} \quad L_xu:= -\ve u_{xx} +a_1u_x \quad \hbox{and} \quad L_yu:= -\ve u_{yy} +a_2u_y;\\
\hbox {and their discrete counterparts} \\
 L^NU(x_i,y_j) =\bigl( \frac{Q_y^C}{\bar k_j} L^N_x+\frac{Q_x^C}{\bar h_i} L^N_y + \bar b_{i,j}\frac{Q^C_xQ_y^C}{\bar h_i\bar k_j}I \bigr) U(x_i,y_j); \\ \\
L^N_xU(x_i,y_j):=\frac{1}{\bar h_i } \bigl (-\ve h_{i+1}D^+_x(\sigma (-\rho _{1:i,j}) D^-_x)+\bar a_1 h_{i}D^-_x\bigr) U(x_i,y_j); \\
 L^N_yU(x_i,y_j):=\frac{1}{\bar k_j} \bigl (-\ve k_{j+1}D^+_y(\sigma (-\rho _{2:i,j}) D^-_y+\bar a_2 k_jD^-_y\bigr) U(x_i,y_j).
\end{eqnarray*}
At each internal mesh point $(x_i,y_j)$, the total truncation error is
\begin{subequations}\label{truncation-error}
\begin{eqnarray}
L^N(U-u) = \frac{Q^C_y}{\bar k_j} \bigl(L_x^N - L_x)u +  \frac{Q^C_x}{\bar h_i} \bigl(L_y^N - L_y)u 
+( \frac{Q^C_y}{\bar k_j}-1) L_xu 
\\
+  (\frac{Q^C_x}{\bar h_i} -1) L_y u 
+\bar b(1- \frac{Q^C_xQ_y^C}{\bar h_i\bar k_j}) u+ \frac{1}{\bar h_i \bar k_j} (\bar f, \psi _{i,j}) - f.
\end{eqnarray}
\end{subequations}
At each internal mesh point $(x_i,y_j)$, the truncation error in one coordinate direction is given by
\begin{subequations}\label{truncation-error-1d}
\begin{eqnarray}
(L_x-L_x^N)u(x_i,y_j)
&=&\ve \frac{(\sigma (-\rho _{1:i+1,j})-1) D_x^+u-(\sigma (-\rho _{1:i,j})-1)D_x^-u}{\bar h_i} \nonumber \\
&+& \ve (\delta _x^2 u - u_{xx} ) +a_1u_x-\bar a_1D_x^-u- \bar a_1\frac{h_i-h_{i+1}}{h_i+h_{i+1}}D_x^-u \nonumber\\
&=&a_1u_x -\bar a_1D_x^-u- \bar a_1\frac{h_i-h_{i+1}}{h_i+h_{i+1}}D_x^-u\\
&-& \ve \frac{(\sigma (-\rho _{1:i,j}) - \sigma (-\rho _{1:i+1,j})) }{\bar h_i} D_x^-u \\
&+& \ve (\delta _x^2 u - u_{xx})  +\ve (\sigma (-\rho _{1:i+1,j})-1)  \delta _x^2 u.
\end{eqnarray}
\end{subequations}
Observe that
\begin{equation}\label{fitting-factor-bounds}
\ve \vert (\sigma (- \rho _{i}) -1) \vert \leq C\min \{ h_i , \ve \};
\quad 
\vert \sigma (-\rho _{i+1})-\sigma ( -\rho _{i})  \vert \leq C\vert \rho _{i+1}- \rho _{i}) \vert .
\end{equation}

\begin{lemma} For all three choices of test functions in (\ref{Fitted-variants-lumped}) we have
\begin{equation}\label{inhom-bound}
\frac{1} {\bar h_i \bar k_j } \vert (\bar f, \psi _{i,j}) - f(x_i,y_j) \vert \leq  C(h_i +k_j+ \frac{h_i-h_{i+1}}{h_i+h_{i+1}}+\frac{k_j-k_{j+1}}{k_j+k_{j+1}}).
\end{equation}
\end{lemma}

\begin{proof} 
For the $\bar L^*$ test functions  (\ref{test_L^*_lumped}),  if $p_i$ denotes the standard hat function
then  we have 
\begin{eqnarray*}
(1, \psi _i - p_i) = h_i \int _{s=0}^1 \frac{1- e^{-\rho _i s}}{1- e^{-\rho _i  }} -s \ ds -  h_{i+1} \int _{s=0}^1 \frac{1- e^{-\rho _{i+1}  s}}{1- e^{-\rho _ {i+1}  }} -s \  ds \\
= h_{i+1}\int _{s=0}^1 \frac{1- e^{-\rho _i  s}}{1- e^{-\rho _i  }} -\frac{1- e^{-\rho _{i+1}  s}}{1- e^{-\rho _{i+1}  }} \ ds +(h_i -  h_{i+1}) \int _{s=0}^1 \frac{1- e^{-\rho _{i}  s}}{1- e^{-\rho _{i}  }} -s \  ds .
\end{eqnarray*} 
We bound these two terms seperately. 
\begin{subequations}\label{star}
\begin{eqnarray}
\Bigl\vert \int _{s=0}^1 \frac{1- e^{-\rho _{i}  s}}{1- e^{-\rho _{i}  }} -s \  ds \Bigr\vert \leq C\min \{ 1 , \rho _{i} \} \\
\hbox{using the inequality} \quad \vert \frac{2(1-x-e^{-x}) +x(1-e^{-x})}{2x(1-e^{-x})} \vert  \leq C \min \{1, x\} . \nonumber \\
\hbox{Also} \quad \Bigl\vert \int _{s=0}^1 \frac{1- e^{-\rho _i  s}}{1- e^{-\rho _i  }} -\frac{1- e^{-\rho _{i+1}  s}}{1- e^{-\rho _{i+1}  }} \ ds  \Bigr\vert \leq C\frac{\vert \bar a_ih_i- \bar a_{i+1}h_{i+1} \vert }{h_{i+1} } \\
\hbox{using } \quad \bigl \vert \frac{d}{dx} \Bigl(  \frac{1}{1-e^{-x}} - \frac{1}{x}   \Bigr) \Bigr\vert  =  \bigl \vert \frac{1}{x^2}  - \frac{e^{-x}}{(1-e^{-x})^2} \Bigr\vert  \leq C \min \{1, \frac{1}{x}\}.\nonumber
\end{eqnarray}
\end{subequations}
Hence, using $h_{i+1} \leq h_i$ and these two inequalities, 
\begin{subequations}
\begin{equation}\label{test-boundx}
 \vert  \frac{ Q^C_x}{\bar h_i} -1 \vert =\vert  \frac{ (1, \psi _i )}{\bar h_i} -1 \vert \leq C( h_{i+1} +\frac{h_i-h_{i+1}}{h_i+h_{i+1}}).
\end{equation}
Likewise,
\begin{equation}\label{test-boundy}
 \vert  \frac{ Q^C_y}{\bar k_j} -1 \vert =\vert  \frac{ (1, \psi ^j )}{\bar k_j} -1 \vert \vert \leq C( k_{j+1} +\frac{k_j-k_{j+1}}{k_j+k_{j+1}}).
\end{equation}
\end{subequations}
Using an analogous argument, one can see that the bounds (\ref{test-boundx}), (\ref{test-boundy}) also apply if the test functions are  
$L$ test functions  (\ref{test_L_lumped}). The bounds (\ref{test-boundx}), (\ref{test-boundy}) are satisfied (trivially) if the test functions are the simple hat functions.

For any choice of the three test functions in (\ref{Fitted-variants-lumped}) 
\begin{eqnarray*}
(\bar f, \psi _{i,j}) - \bar h_i \bar k_j f(x_i,y_j)) = (\bar f-f(x_i,y_j), \psi _{i,j} ) -f(x_i,y_j) (1,\psi _{i,j} -p _{i,j}) \\
\vert  (\bar f-f(x_i,y_j), \psi _{i,j} ) \vert \leq C(h_i+k_j) \bar h_i \bar k_j
\\
(1,\psi _{i,j} -p_{i,j})  = (1,\psi _{i} -p_{i})(1, \psi ^j) + (1,p _{i})(1, \psi ^j-p ^{j}))\\
 \leq C(h_i +k_j+ \frac{h_i-h_{i+1}}{h_i+h_{i+1}}+\frac{k_j-k_{j+1}}{k_j+k_{j+1}})\bar h_i \bar k_j.
\end{eqnarray*}
\end{proof}

The discrete solution $U$ can be decomposed in an analogous fashion to the continuous solution. We write
\begin{eqnarray*}
U=V+W_R+W_T+W_{RT}, \quad \hbox{where the regular component satisfies } \\
L^NV(x_i,y_j)=f(x_i,y_j), \ (x_i,y_j) \in \Omega ^N , \quad V(x_i,y_j) =v(x_i,y_j)\ (x_i,y_j) \in \partial \Omega ^N ;  \\
\hbox{and each of the layer functions satisfy }  \\ L^NW(x_i,y_j) =0, \ (x_i,y_j) \in \Omega ^N , \quad W(x_i,y_j) =w(x_i,y_j)\ (x_i,y_j) \in \partial  \Omega ^N.
\end{eqnarray*}

\begin{lemma} For the regular component, the nodal error bound 
\[
\vert (V - v)(x_i,y_j) \vert \leq C N^{-1}, \quad (x_i,y_j) \in \Omega ^N
\]
is satisfied for all three schemes in (\ref{Fitted-variants-lumped}).
\end{lemma}
\begin{proof}
We first bound the truncation error for the regular component and then use the discrete maximum principle and a suitable discrete barrier function to deduce the error bound. 

Combining the bounds in (\ref{test-boundx}), (\ref{test-boundy}), (\ref{inhom-bound}) into the expression for the truncation error (\ref{truncation-error}), we deduce that
\begin{eqnarray*}
\vert L^N(V-v)(x_i,y_j) \vert \leq C  \vert (L_x^N - L_x)v (x_i,y_j)\vert +  C \vert (L_y^N - L_y)v (x_i,y_j)\vert +
\\
C(1+\Vert L_xv \Vert+\Vert L_y v \Vert ) (h_i +k_j+ \frac{h_i-h_{i+1}}{h_i+h_{i+1}}+\frac{k_j-k_{j+1}}{k_j+k_{j+1}}).
\end{eqnarray*}
Using the bounds (\ref{fitting-factor-bounds}) in the one dimensional truncation error (\ref{truncation-error-1d}), we also have that
\[
\vert (L_x-L_x^N)v(x_i,y_j) \vert  \leq C(1+\Vert v_x \Vert + \Vert v_{xx} \Vert)  (h_i+ \frac{h_i-h_{i+1}}{h_i+h_{i+1}}) +Ch_i\ve \Vert v_{xxx} \Vert.
\]
Hence, using the bounds (\ref{regular-bnd})  on the derivatives of the regular component $v$,  we have establshed the truncation error bound
\[
\vert L^N(V-v) (x_i,y_j) \vert \leq C(h_i +k_j+ \frac{h_i-h_{i+1}}{h_i+h_{i+1}}+\frac{k_j-k_{j+1}}{k_j+k_{j+1}}).
\]
Observe that the finite difference operator $L^N_x$ can be rewritten in the form
\[
L^N_xU(x_i,y_j)=\frac{\ve}{\bar h_i } \bigl (\sigma (\rho _{1:i,j}) D^-_x-\sigma (-\rho _{1:i+1,j}) D^+_x \bigr) U(x_i,y_j).
\]
Consider the one dimensional barrier function
\[
T_x(x_i) = e^{-\alpha\frac{1-\sigma _x-x_i}{\ve}},\ x _i \leq 1-\sigma_x ; \qquad T_x(x_i)=1, x _i \geq 1 -\sigma_x. 
\]
In the coarse mesh area where $ x_i < 1-\sigma _x$: 
\begin{eqnarray*}
L^N_xT_x(x_i) \geq \frac{\ve}{H } \bigl (\sigma (\rho ) D^-_x-\sigma (-\rho ) D^+_x \bigr) T(x_i) = 0 ,  \ \rho := \frac{\alpha H}{\ve}; \quad \hbox{as} \quad \sigma '(x) >0.
\end{eqnarray*}
At the transition point
\[
 L_x^N T_x(1-\sigma_x) \geq \frac{\alpha }{\bar h_i} \geq CN
\]
and $L_x^N T_x(x_i) =0$ in the fine mesh area where $x_i>1-\sigma_x$.
We also have that, for $N$ sufficiently large, 
\begin{eqnarray*}
L^N_x(x_i) = \frac{\ve}{\bar h_i} (\sigma_x (\rho_{1:i,j})) - \sigma _x (-\rho_{1:i+1,j}))  \geq 0, \quad \forall i; \\
L^N_x(x_i) = \frac{\ve}{\bar h_i} (\rho_{1:i,j}+ (\sigma_x (-\rho_{1:i,j})) - \sigma _x (-\rho_{1:i+1,j}))\geq \frac{a_1(x_i,y_j)}{2} \quad i \neq N/2.
\end{eqnarray*}
Then, combining these two functions $x_i, T_x(x_i)$ into one barrier function, we deduce that
\[
\vert (V - v)(x_i,y_j)  \vert \leq CN^{-1}(x_i+T_x(x_i)+ y_j+T_y(y_j))\leq  CN^{-1}.
\]
\end{proof}
Consider the following fitted finite difference operator
\[
L^N_{x,1}:= -\ve \sigma  (-\rho ^*_j) \delta _x^2 +a(1,y_j) D_x^-, \quad \rho ^*_j := \frac{a(1,y_j)h}{\ve}.
\]
On the fine mesh region $(1-\sigma_x,1)\times (0,1) $
\begin{equation}\label{key}
L^N_{x,1} E^* (x_i,y_j)  =0, \quad \hbox{where} \quad E ^* (x_i,y_j) := e^{-\frac{a(1,y_j) (1-x_i)}{\ve}},
\end{equation}
which is a key property of the  fitted operator $L^N_x$. 
\begin{lemma} For the layer components, the error bounds
\begin{subequations}
\begin{eqnarray} 
\vert (W_R-w_R)(x_i,y_j) \vert \leq C N^{-1},\quad (x_i,y_j) \in \Omega ^N;\\ 
\vert (W_T-w_T)(x_i,y_j) \vert \leq C N^{-1},\quad (x_i,y_j) \in \Omega ^N;\\ 
\vert (W_{TR}-w_{TR})(x_i,y_j) \vert \leq C N^{-1},\quad (x_i,y_j) \in \Omega ^N;
\end{eqnarray}
are satisfied by all three schemes in (\ref{Fitted-variants-lumped}).
\end{subequations}

\end{lemma}

\begin{proof}
Consider the barrier function
\[
B_x(x):= e^{-\frac{\alpha _1 (1-x)}{2\ve}}, \quad \rho := \frac{\alpha _1h}{\ve},\  h: = \min h_i.
\]
Then at all internal mesh points
\begin{eqnarray*}
L_x^N B_x(x_i)  
&\geq&   \frac{\alpha_1}{2 \bar h_i} \Bigl (\frac{\sigma (\rho _i)}{\sigma (\rho/2 )}  - \frac{\sigma (-\rho _{i+1})}{\sigma (-\rho/2 )} \Bigr)  B_x(x_i) \\
&\geq& CN(1-e^{-\rho/2}) e^{-\frac{\alpha _1 (1-x_i)}{2\ve}} > 0, \quad \hbox{as}\quad  \sigma'(x) >0.
\end{eqnarray*}
Hence we have the following bound on the discrete singular term $W_R$
\[
\vert W_R(x_i,y_j) \vert  \leq C e^{-\frac{\alpha _1(1-x_i)}{2\ve}}.
\]
Observe that in the case where the horizontal mesh is a uniform mesh (i.e., $\sigma _x =0.5$), then $ e^{-\frac{\alpha _1\sigma _x}{2\ve}}\leq  N^{-1}$.

Also, by the choice of the transition parameter $\sigma _x$ in (\ref{transition parameters}) and the pointwise bound (\ref{edge-bnd}) on the layer component
\[
\vert w_R(x_i,y_j) \vert \leq CN^{-2},\quad  x_i \leq 1-\sigma _x.
\] 
Hence,  for the mesh points outside the right boundary layer region,  
\[
\vert (W_R-w_R)(x_i,y_j) \vert \leq C N^{-1},\quad  x_i \leq 1-\sigma_x.
\]
For mesh points within the side region $(1-\sigma _x,1) \times (0,1)$, we have that the truncation error, along each level $y=y_j$, is
\begin{eqnarray*}
(L_x-L_x^N)w_R
&=&(a_1(x_i,y_j)-a_1(1,y_j)) \Bigl( \frac{ \partial  w_{R}}{\partial x}-D_x^-w_R\Bigr)  +(a_1-\bar a_1)D_x^-w_R\\
&-& \ve \frac{(\sigma (-\rho _{1:i,j}) - \sigma (-\rho _{1:i+1,j})) }{ h} D_x^-w_R + \ve \bigl(\sigma (-\rho _{1:i+1,j})- \sigma (-\rho _j^*)\bigr) \delta _x^2 w_R  \\
&+& \ve \bigl(\sigma (-\rho ^* _{j})\delta _x^2 w_R - \frac{ \partial ^2 w_{R}}{\partial x^2} \bigr) +a_1(1,y_j) \frac{ \partial  w_{R}}{\partial x} - a_1(1,y_j))D_x^-w_R.
\end{eqnarray*}
We next bound some of the terms in this expression
\begin{eqnarray*}
\Bigl \vert (a_1-\bar a_1)D_x^-w_R(x_i,y_j)\Bigr \vert &\leq & C\Bigl \vert  \int _{x_{i-1}}^{x_i} \frac{ \partial  w_{R}}{\partial s}(s,y_j) \ ds \Bigr \vert  \\
&\leq & Ce^{-\frac{\alpha _1(1-x_i)}{\ve}}(1- e^{-\frac{\alpha _1 h }{\ve}}). \\
\Bigl \vert (a_1(x_i,y_j)-a_1(1,y_j)) \Bigl( \frac{ \partial  w_{R}}{\partial x}-D_x^-w_R(x_i,y_j)\Bigr) \Bigr \vert &\leq & C\frac{(1-x_i) }{\ve}  e^{-\frac{\alpha _1(1-x_i)}{\ve}}(1- e^{-\frac{\alpha _1h }{\ve}}) \\
&\leq & Ce^{-\frac{\alpha _1(1-x_i)}{2\ve}}(1- e^{-\frac{\alpha _1 h }{\ve}}).\\
\Bigl \vert  \ve \frac{(\sigma (-\rho _{1:i,j}) - \sigma (-\rho _{1:i+1,j})) }{ h} D_x^-w_R (x_i,y_j)  \Bigr \vert &\leq &Ce^{-\frac{\alpha _1(1-x_i)}{\ve}}(1- e^{-\frac{\alpha _1h }{\ve}}). \\
\Bigl \vert \ve \bigl(\sigma (-\rho _{1:i+1,j})- \sigma (-\rho _j^*)\bigr) \delta _x^2 w_R (x_i,y_j)\Bigr \vert & \leq & C\frac{h(1-x_i) }{\ve^2}  e^{-\frac{\alpha _1 (1-x_i)}{\ve}} \\ &\leq & Ce^{-\frac{\alpha _1(1-x_i)}{2\ve}}(1- e^{-\frac{\alpha_1 h }{\ve}}). 
\end{eqnarray*}
The following terms remain to be bounded
\[
T:= \ve \sigma (-\rho ^* _{j})\delta _x^2 w_R- a_1(1,y_j)D_x^-w_R + \bigl(  - \ve \frac{ \partial ^2 w_{R}}{\partial x^2} + a_1(1,y_j) \frac{ \partial  w_{R}}{\partial x} \bigr).
\]
Based on the arguments in \cite[pg. 1765]{hemkera}, we can decompose the layer component as follows
\[
w_R(x,y) = w_R(1,y)  e^{-\frac{a_1(1,y) (1-x)}{\ve}} +\ve z_R(x,y);
\]
where (using  the inequality $t^ne^{-t} \leq C e^{-\theta  t}, n \geq 1, t \geq 0, 0 < \theta <1$ instead of using the inequality $t^ne^{-t} \leq C e^{-t/2}, n \geq 1, t \geq 0$ and the localized derivative bounds in \cite[pp.132-134]{lady})
we can deduce that
\begin{equation}\label{orth-bnd-remainder}
\Bigl \vert \frac{\partial ^{i+j}z_R}{ \partial  x^i\partial  y^j} (x,y) \Bigr\vert  \leq
C\ve ^{-i}(1+\ve ^{-j})e^{-\frac{\alpha _1(1-x)}{2\ve}}, \ i,j \leq 3.
\end{equation}
From (\ref{key}), we have
\[
T = \ve \Bigl(  \ve \sigma (-\rho ^* _{j})\delta _x^2 z_R- a_1(1,y_j)D_x^-z_R + \bigl(  - \ve \frac{ \partial ^2 z_{R}}{\partial x^2} + a_1(1,y_j) \frac{ \partial  z_{R}}{\partial x} \bigr)\Bigr)
\]
and 
\begin{eqnarray*}
\Bigl \vert \Bigl( \frac{ \partial  z_{R}}{\partial x}-D_x^-z_R\Bigr) \Bigr \vert &\leq & C\frac{1}{\ve} e^{-\frac{\alpha _1(1-x_i)}{2\ve}}(1- e^{-\frac{\alpha _1h }{\ve}});\\
\Bigl \vert  \ve \frac{(\sigma (-\rho _{1:i,j}) - \sigma (-\rho _{1:i+1,j})) }{ h} D_x^-z_R   \Bigr \vert &\leq &Ce^{-\frac{\alpha _1(1-x_i)}{2\ve}}(1- e^{-\frac{\alpha _1h }{\ve}}) ;\\
\Bigl \vert \ve \bigl( \sigma (-\rho _j^*)-1\bigr) \delta _x^2 z_R \Bigr \vert & \leq & C\frac{h}{\ve^2}  e^{-\frac{\alpha _1 (1-x_i)}{2\ve}} \leq C\frac{1}{\ve}e^{-\frac{\alpha _1(1-x_i)}{2\ve}}(1- e^{-\frac{\alpha_1 h }{\ve}}); \\
\Bigl \vert \ve \bigl(\frac{ \partial ^2 z_{R}}{\partial x^2} - \delta _x^2 z_R \bigr)\Bigr \vert & \leq &  C\frac{1}{\ve}e^{-\frac{\alpha _1(1-x_i)}{2\ve}}(1- e^{-\frac{\alpha _1h }{\ve}}) . 
\end{eqnarray*}
Collecting all of the bounds above, we have the following truncation error bound
\[
\vert (L_x-L_x^N)w_R(x_i,y_j) \vert \leq  Ce^{-\frac{\alpha _1 (1-x_i)}{2\ve}}(1- e^{-\frac{\alpha _1 h }{\ve}}), \quad x_i > 1-\sigma _x.
\]
Furthermore, using (\ref{orth-bnd}), we have that 
\[
\vert (L_y-L_y^N)w_R(x_i,y_j) \vert \leq  C\frac{N^{-1}}{\ve} e^{-\frac{\alpha _1(1-x_i)}{\ve}}.
\]
Hence, we arrive at the truncation error bound:
\[
\vert L^N(W_R-w_R) (x_i,y_j) \vert \leq C\Bigl( \frac{N^{-1}}{\ve} +(1- e^{-\frac{\alpha h }{\ve}}) \Bigr)e^{-\frac{\alpha (1-x_i)}{2\ve}}, \quad x_i > 1-\sigma _x.
\]
In the fine mesh region where $x_i > 1-\sigma _x$, we have that
\begin{eqnarray*}
L^N B_x(x_i)  \geq C\frac{1}{h}(1-e^{-\rho _1/2}) e^{-\frac{\alpha _1(1-x_i)}{2\ve}}, \qquad B_x(x):= e^{-\frac{\alpha _1(1-x)}{2\ve}}.
\end{eqnarray*}
Thus
\[
\vert (W_R-w_R) (x_i,y_j) \vert \leq C\Bigl( \frac{N^{-1} \rho }{(1- e^{-\frac{\alpha _1h }{\ve}})} +N^{-1} \Bigr)B_x(x_i) \leq CN^{-1}.
\]
The bound on $\vert (W_T-w_T) (x_i,y_j) \vert$ is also deduced using an analogous argument. 

For the corner layer function, we can establish that
\[
w_{RT}(x,y) = w_R(1,1)  e^{-\frac{a_1(1,1) (1-x)}{\ve}} e^{-\frac{a_2(1,1) (1-y)}{\ve}} +\ve z_{RT}(x,y);
\]
where
\begin{equation}\label{orth-bnd-remainder-corner}
\Bigl \vert \frac{\partial ^{i+j}z_{RT}}{ \partial  x^i\partial  y^j} (x,y) \Bigr\vert  \leq
C\ve ^{-(i+j)}e^{-\frac{\alpha _1(1-x)}{2\ve}}e^{-\frac{\alpha _2(1-x)}{2\ve}}, \ i,j \leq 3.
\end{equation}
Then repeating the argument from above, we deduce that
\[
\vert (W_{RT}-w_{RT}) (x_i,y_j) \vert  \leq CN^{-1}.
\]
\end{proof}

By combining the results from the previous two lemmas, we arrive at the nodal error bound
\begin{theorem} {\it (Nodal convergence)}
If $u$ is the solution of (\ref{prob-2d-elliptic}), then
\begin{equation}
\vert (U-u)(x_i,y_j) \vert \leq C N^{-1},\quad (x_i,y_j) \in \Omega ^N;
\end{equation}
where $U$ is the numerical approximation generated by any one of the  three schemes in (\ref{Fitted-variants-lumped}).
\end{theorem}

 On the Shishkin mesh, this is easily extended to a global error bound using simple bilinear interpolation. 
If $p_i(x) (p^j(y))$ denotes the standard hat function centered at $x=x_i$ ($y=y_j$) and
we denote the bilinear interpolants of the exact solution $u$ and the numerical solution $U$ by $u_{I, BL}$ and $U_{I, BL}$; then
\[
u_{I, BL}(x,y) := \sum _{i,j=1}^{N} u(x_i,y_j) p_i(x) p^j(y), \ U_{I, BL}(x,y) := \sum _{i,j=1}^{N} U(x_i,y_j) p_i(x) p^j(y).
\]
Note that if the trial functions are chosen as bilinear  functions, then $u_I \equiv u_{I, BL}$ and $U \equiv U_{I, BL}$.
Using the triangle inequality and the  interpolation bound \cite[Theorem4.2]{styor4} 
\[
\Vert u -u_{I, BL}\Vert _\infty \leq C (N^{-1} \ln N)^2
\]
on the Shishkin mesh, we easily deduce the following global error bound.
\begin{theorem}\label{main-result} {\it (Global convergence)}
If $u$ is the solution of (\ref{prob-2d-elliptic}), then
\begin{equation}
\Vert U-u \Vert _\infty \leq C N^{-1},
\end{equation}
where $U$ is the numerical approximation generated by either the  numerical method (\ref{test_L^*_lumped}) or the  numerical method  (\ref{test_L_lumped}).
\end{theorem}
For notational simplicity, we have taken the same number $N$ of elements in each coordinate direction. If one uses $N$ mesh elements in the horizontal direction and $M$ elements in the vertical direction, 
then one can easily establish $\Vert U-u \Vert _\infty \leq C N^{-1} +CM^{-1}.$

\begin{remark}
The numerical method (\ref{test_hat_lumped}) uses exponential basis functions in the trial space. Under additional regularity assumptions (assume $u \in C^{4,\gamma }(\overline{ \Omega })$ and that the pointwise bounds \cite[(2.8d)-(2.8f)]{dervs-elliptic} on the layer components $w_R,w_T, w_{TR}$ are valid), then one can apply the arguments in \cite[Lemma 5.1]{orst4} to each of the components $v, w_R,w_T, w_{TR}$ separately to deduce that $\Vert u-u_I \Vert _\infty \leq C N^{-1}$ for $\bar L$-splines on the Shishkin mesh. However, for practical reasons in our evaluation of the global accuracy of the numerical approximations, we confine our attention to bilinear interpolants in the numerical section of this paper. Note that for method (\ref{test_hat_lumped}) the bilinear interpolant $U_{I, BL}$ will satisfy the global error bound in Theorem \ref{main-result}.
\end{remark}

In the non-singularly perturbed case, we  see that the schemes return to classical Galerkin with bilinear elements.
Hence, we can establish the following result.
\begin{theorem} Assume that $u \in C^{4,\gamma} (\bar \Omega)$. 
If $\ve \geq \ve _0 >0$, where $\ve _0$ is a fixed constant, then for all three fitted schemes in (\ref{Fitted-variants-lumped}) 
\begin{equation}
\Vert  U-u\Vert _\infty\leq C N^{-2},
\end{equation}
if $N$ is sufficiently large. 
\end{theorem}

\begin{proof}When $\ve \geq \ve _0$, then 
\[
\Bigl \Vert \frac{\partial ^{i+j}u}{ \partial  x^i\partial  y^j}  \Bigr\Vert  \leq
C\ve _0 ^{-(i+j)}\leq C, \quad i,j \leq 4
\] 
and, then, $
\Vert u -u_I \Vert _\infty \leq CN^{-2}.
$

Let us now examine the nodal error. If $\ve \geq \ve _0 >0$ and $N$ is sufficiently large, then the mesh is uniform ($h_i=h,\ \forall i; k_j=k, \ \forall j$).

We first establish the result for the scheme (\ref{test_hat_lumped}), which  simplifies 
\begin{eqnarray*}
\bigl(-\ve ( D^+_x(\sigma (-\rho _{1:i,j}) D^-_x)+ D^+_y(\sigma (-\rho _{2:i,j}) D^-_y) +\bar a_1 D^-_x
 +\bar a_2 D^-_y\bigr) U(x_i,y_j)\\
+ b(x_i,y_j)   U(x_i,y_j)=   \sum _{n=i-1}^{i+1} \sum _{m=j-1}^{j+1} \gamma _{n,m} f(x_n,y_m) ;
\\
\hbox{where} \qquad 
\gamma _{i-1,j-1} =\gamma _{i+1,j-1} =\gamma _{i-1,j+1} =\gamma _{i+1,j+1} =  \frac{1}{16};\\
\gamma _{i,j-1} =\gamma _{i,j-1} =\gamma _{i-1,j} =\gamma _{i+1,j} =  \frac{1}{4}; \quad \gamma _{i,j} =  \frac{1}{2}.
 \end{eqnarray*}
A standard Taylor series expansion yields
\begin{equation}\label{classical-RHS}
\bigl\vert \sum _{n=i-1}^{i+1} \sum _{m=j-1}^{j+1} \gamma _{n,m} f(x_n,y_m)  -f(x_i,y_j)  \bigr\vert \leq CN^{-2}.
\end{equation}
Note also that, on a uniform mesh, 
\[
D^-_x= D_x^0 - \frac{h}{2} \delta ^2_x; \quad \hbox{where} \quad
\quad
D^0_xU(x_i,y_j) := \frac{U(x_i,y_j) - U(x_i,y_j)}{2h}.
\]
Also
\begin{eqnarray*}
-\ve  D^+_x(\sigma (-\rho _{1:i,j}) D^-_x)+  \bar a_{1:i,j} D^-_x = -\ve  \delta^2_x + a_{1}(x_i,y_j) D^0_x\\
+  (\bar a_{1:i,j} -a_1(x_i,y_j))D^0_x+ \frac{\bar a_{1:i+1,j} -\bar a_{1:i,j}}{2} D^+_x + \\
\frac{\ve}{h} \Bigl( (\sigma (-\rho _{1:i,j}) +\frac{\rho _{1:i,j}}{2} -1)D^-_x) -  (\sigma (-\rho _{1:i+1,j}) +\frac{\rho _{1:i+1,j}}{2} -1)D^+_x)\Bigr);
 \end{eqnarray*}
and
\begin{eqnarray*}
 (\bar a_{1:i,j} -a_1(x_i,y_j))D^0_x u(x_i,y_j) =- \frac{h}{2} \frac{\partial a_1}{\partial x}(x_i,y_j) \frac{\partial u}{\partial x} (x_i,y_j)+O(h^2);\\
 \frac{\bar a_{1:i+1,j} -\bar a_{1:i,j}}{2} D^+_x u(x_i,y_j) = \frac{h}{2}
 \frac{\partial a_1}{\partial x}(x_i,y_j) \frac{\partial u}{\partial x} (x_i,y_j) +O(h^2);\\
\sigma (-\rho _{1:i,j}) +\frac{\rho _{1:i,j}}{2} -1 =\frac{\rho _{1:i,j}}{2}\coth \Bigl(\frac{\rho _{1:i,j}}{2} \Bigr)  -1 = \frac{\rho ^2_{1:i,j}}{12} +O(h^4);\\
\frac{\ve \rho ^2_{1:i,j}}{h} (D^+_x-D^-_x)u(x_i,y_j) =\ve \rho ^2_{1:i,j}\delta ^2_xu(x_i,y_j) = O(h^2).
\end{eqnarray*}
Collecting all these terms, we can deduce the truncation error bound
\[
\vert L^N(u-U) (x_i,y_j) \vert \leq Ch^2 +Ck^2,
\]
and the nodal error bound follows, using the discrete maximum principle. This completes the proof in the case of the scheme (\ref{test_hat_lumped})

In the case of the other two fitted schemes, use the inequality
\[
\Bigl \vert \frac{1-e^{-sx}}{1-e^{-x} } - s \Bigr \vert \leq Cx^2, \quad 0\leq s \leq 1, \ x >0
\]
to establish that, when $\ve \geq \ve _0$, that
\[
 \vert  \frac{ Q^C_x}{h} -1 \vert \leq Ch^2 \quad \hbox{and} \quad   \vert  \frac{ Q^C_y}{ k} -1 \vert \leq Ck^2.
\]
Complete the proof as above. 

\end{proof}

 \section{Numerical examples}

In this final section, we will estimate the global (as opposed to the nodal) accuracy of our numerical methods.

{\bf Example 1}
Consider the following constant coefficient test problem:
\begin{equation}\label{test1}
-\ve \triangle u + 2u_x+3u_y = f,  (x,y) \in \Omega; \ u =0,\ (x,y) \in \partial \Omega ,                 
\end{equation}
where $f$ is such that the exact solution is 
\[
u(x,y) = \Bigl( x - \frac{e^{2x/\varepsilon} -1}{e^{2/\varepsilon} -1}\bigr) \bigl(y - \frac{e^{3y/\varepsilon} -1}{e^{3/\varepsilon} -1} \bigr)+\sin (\pi x) \sin (\pi y).
\]
A sample plot of the computed solution using  the numerical scheme (\ref{test_L^*_lumped}) is displayed in Figure 1 and the corresponding global error is displayed in  Figure 2. 
\begin{figure}[ht!]
\centering
\includegraphics[width=0.6\textwidth]{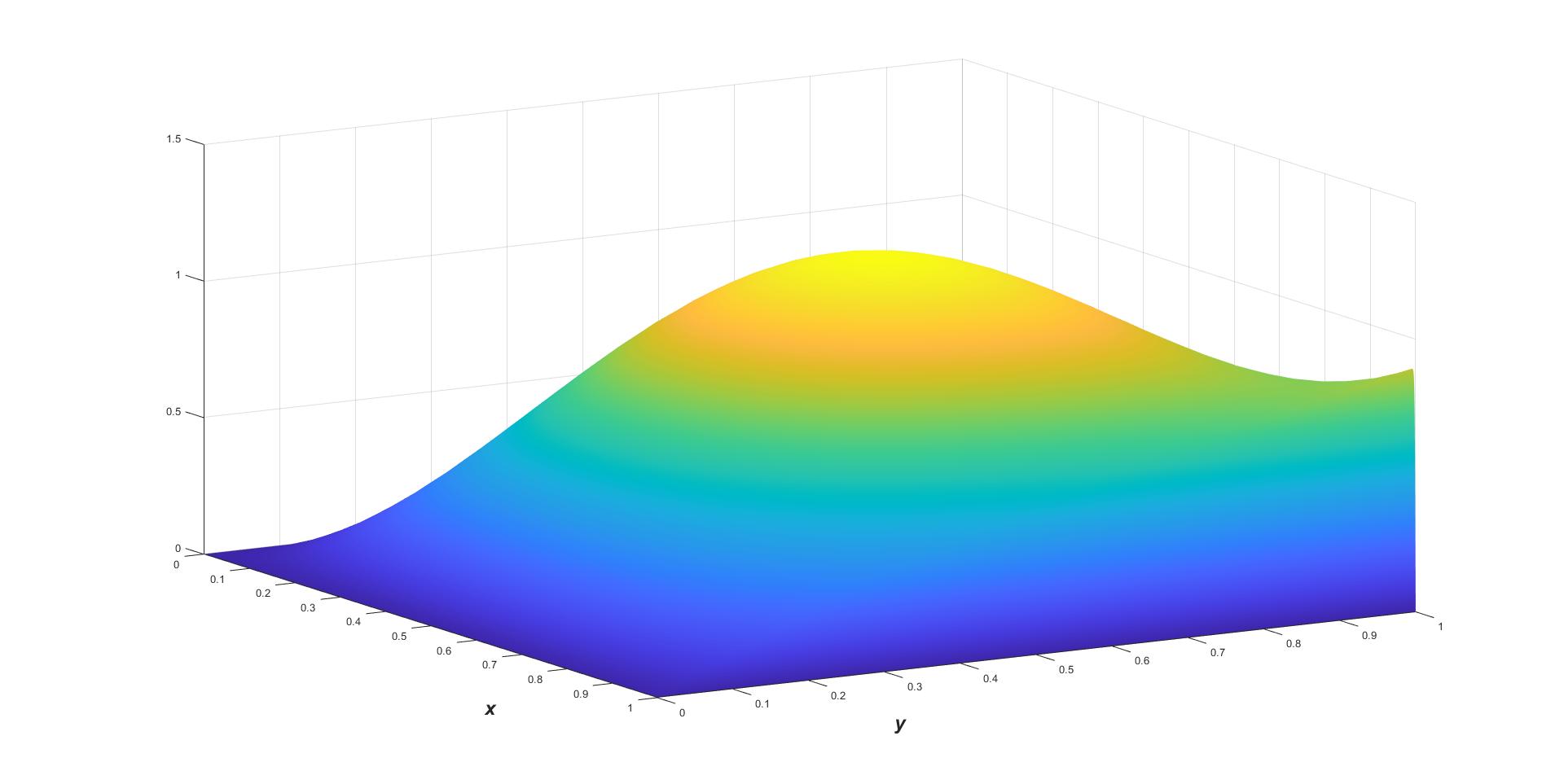}%
\caption{Computed solution with the numerical scheme (\ref{test_L^*_lumped}) applied to problem (\ref{test1})  for $\ve=2^{-10}$ and $N=M=64$}%
\label{figure1}%
\end{figure}

 \begin{figure}[ht!]
\centering
\includegraphics[width=0.6\textwidth]{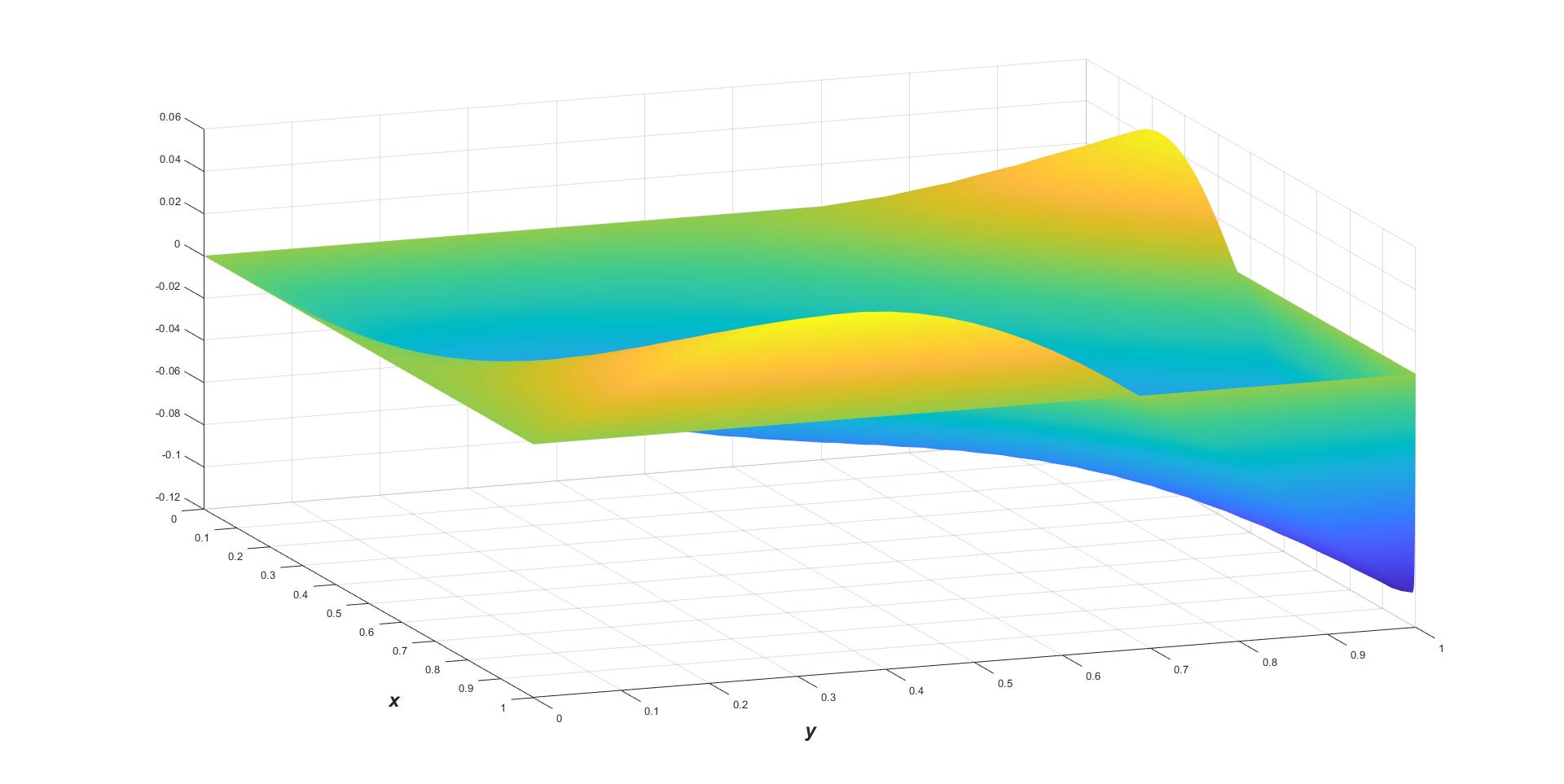}%
\caption{Global error with the numerical scheme (\ref{test_L^*_lumped})  applied to problem (\ref{test1}) for $\ve=2^{-10}$ and $N=M=64$}%
\label{figure2}%
\end{figure}

In the case of the three fitted schemes and simple upwinding on the Shishkin mesh, the global errors  are estimated by calculating the maximum error over a fine Shishkin mesh. That is, 
\[
\Vert U -u \Vert \approx E^N_\ve := \max _{(x_i,y_j) \in \Omega ^{2048}_S} \vert (U_{I,BL}-u)(x_i,y_j) \vert
\]
The  global errors in Tables  1,2,3 for the fitted schemes  (\ref{test_L^*_lumped}),   (\ref{test_L_lumped})  and (\ref{test_hat_lumped}) indicate first order convergence for this constant coefficient problem
and second order convergence for $\ve$ of order one. These errors can be compared to the global errors in Table 4 for simple upwinding on the same Shishkin mesh. For large values of $\ve$, the fitted schemes are a significant improvement on basic upwinding. Overall, the fitted scheme  (\ref{test_L^*_lumped}) performs best for this test problem. These numerical results are in agreement with the theoretical error bounds established in Theorem 3.   We  estimate the global  orders of local convergence using the double-mesh principle \cite[\S 8.6]{fhmos}.
The orders of global convergence for (\ref{test_L^*_lumped}) in Table 5 can be compared to the global orders of convergence of simple upwinding in Table 6.

{\bf Example 2} 
Consider the variable coefficient problem
\begin{equation}\label{test3}
-\ve \triangle u + (2+x+x^2+3xy)u_x+(3+y+y^2+2xy)u_y= 16x(1-x)y(1-y))                   
\end{equation}
with $u=0$ on the boundary. The orders of global convergence in Table 7  for the fitted scheme   (\ref{test_L^*_lumped}) indicate first order convergence. These orders should be compared to the orders for standard upwinding  in Table 8  on the same Shishkin mesh, where the $\log $ defect in the orders is  evident. 

\begin{table}\label{T1}
\caption{Global errors $E^N_\ve $ for the scheme (\ref{test_L^*_lumped}) applied to problem  (\ref{test1})} 
\begin{tabular}{|r|rrrrrrrr|}\hline 
$\varepsilon| N$&$N=8$& 16&32&64&128&256&512&1024\\ \hline 
$2^{0}$&   0.0666 &   0.0170 &   0.0043  &   0.0011  &   0.0003 &   0.0001 &   0.0000 &   0.0000\\ 
$2^{-4}$ &   0.3201 &   0.1305 &   0.0408  &   0.0098  &   0.0022 &   0.0005 &   0.0001 &   0.0000\\ 
$2^{-8}$  &   0.4463 &   0.3148 &   0.1831  &   0.0957  &   0.0444 &   0.0171 &   0.0052 &   0.0014\\ 
$2^{-12}$  &   0.4601 &   0.3296 &   0.1982  &   0.1103  &   0.0581 &   0.0295 &   0.0145 &   0.0069\\ 
$2^{-16}$&   0.4610 &   0.3305 &   0.1992  &   0.1112  &   0.0590 &   0.0303 &   0.0153 &   0.0077\\ 
$2^{-20}$ &   0.4610 &   0.3306 &   0.1993  &   0.1113  &   0.0591 &   0.0304 &   0.0154 &   0.0077\\ \hline 
\end{tabular}
\end{table}

\begin{table}\label{T2}
\caption{Global errors $E^N_\ve $ for the scheme (\ref{test_L_lumped}) applied to  problem  (\ref{test1})} 
\begin{tabular}{|r|rrrrrrrr|}\hline 
$\varepsilon| N$&$N=8$& 16&32&64&128&256&512&1024\\ \hline 
$2^{0}$ &   0.0754 &   0.0196 &   0.0049  &   0.0012  &   0.0003 &   0.0001 &   0.0000 &   0.0000\\ 
$2^{-4}$&   0.9269 &   0.3724 &   0.1075  &   0.0260  &   0.0059 &   0.0013 &   0.0003 &   0.0001\\ 
$2^{-8}$ &   1.1373 &   0.7007 &   0.3563  &   0.1703  &   0.0771 &   0.0303 &   0.0096 &   0.0026\\ 
$2^{-12}$  &   1.1545 &   0.7246 &   0.3769  &   0.1911  &   0.0980 &   0.0493 &   0.0243 &   0.0116\\ 
$2^{-16}$ &   1.1556 &   0.7261 &   0.3782  &   0.1926  &   0.0993 &   0.0506 &   0.0255 &   0.0128\\ 
$2^{-20}$ &   1.1557 &   0.7262 &   0.3783  &   0.1927  &   0.0994 &   0.0506 &   0.0256 &   0.0128\\ \hline 
\end{tabular}
\end{table}
\begin{table}{ \label{T3}
\caption{Global errors $E^N_\ve $ for the scheme (\ref{test_hat_lumped}) applied to  problem  (\ref{test1})} 
\begin{tabular}{|r|rrrrrrrr|}\hline 
$\varepsilon| N$&$N=8$& 16&32&64&128&256&512&1024\\ \hline 
$2^{0}$ &   0.0710 &   0.0181 &   0.0046  &   0.0011  &   0.0003 &   0.0001 &   0.0000 &   0.0000\\ 
$2^{-4}$ &   0.5840 &   0.2380 &   0.0703  &   0.0173  &   0.0039 &   0.0009 &   0.0002 &   0.0000\\ 
$2^{-8}$  &   0.7150 &   0.4598 &   0.2575  &   0.1302  &   0.0595 &   0.0231 &   0.0072 &   0.0019\\ 
$2^{-12}$ &   0.7200 &   0.4798 &   0.2768  &   0.1479  &   0.0759 &   0.0381 &   0.0187 &   0.0089\\ 
$2^{-16}$ &   0.7202 &   0.4811 &   0.2781  &   0.1490  &   0.0770 &   0.0391 &   0.0196 &   0.0098\\ 
$2^{-20}$ &   0.7203 &   0.4812 &   0.2781  &   0.1491  &   0.0770 &   0.0391 &   0.0197 &   0.0099\\ \hline 
\end{tabular}}
\end{table}
\begin{table}\label{T4}
\caption{Global errors $E^N_\ve $ for simple upwinding applied to  problem  (\ref{test1})} 
\begin{tabular}{|r|rrrrrrrr|}\hline 
$\varepsilon| N$&$N=8$& 16&32&64&128&256&512&1024\\ \hline 
$2^0$ &   0.1245 &   0.0703 &   0.0372  &   0.0191  &   0.0097 &   0.0049 &   0.0024 &   0.0012\\ 
$2^{-4}$ &   0.6516 &   0.3804 &   0.2069  &   0.1116  &   0.0598 &   0.0319 &   0.0170 &   0.0090\\ 
$2^{-8}$ &   0.7913 &   0.4592 &   0.2569  &   0.1417  &   0.0780 &   0.0422 &   0.0226 &   0.0120\\ 
$2^{-12}$ &   0.7994 &   0.4653 &   0.2602  &   0.1434  &   0.0790 &   0.0428 &   0.0229 &   0.0122\\ 
$2^{-16}$ &   0.7999 &   0.4656 &   0.2604  &   0.1436  &   0.0790 &   0.0428 &   0.0229 &   0.0122\\ 
$2^{-20}$ &   0.7999 &   0.4657 &   0.2604  &   0.1436  &   0.0790 &   0.0428 &   0.0229 &   0.0122\\ \hline 
\end{tabular}
\end{table}

\begin{table}\label{T5}
\caption{Orders of global convergence for the scheme (\ref{test_L^*_lumped})  applied to  problem  (\ref{test1})} 
\begin{tabular}{|r|rrrrrrr|}\hline 
$\varepsilon| N$&$N=8$& 16&32&64&128&256&512\\ \hline 
$2^{0}$ &   1.9025 &   1.9762 &   1.9943  &   1.9987  &   1.9997 &   2.0000 &   2.0000\\ 
$2^{-4}$  &   1.2654 &   1.6719 &   1.9879  &   1.7876  &   1.6014 &   1.6438 &   1.6834\\ 
$2^{-8}$  &   0.8021 &   0.5796 &   0.8045  &   0.9202  &   1.2133 &   1.6216 &   1.8799\\ 
$2^{-12}$ &   0.7724 &   0.5813 &   0.7981  &   0.8801  &   0.9403 &   0.9742 &   0.9906\\ 
$2^{-16}$ &   0.7704 &   0.5813 &   0.7975  &   0.8797  &   0.9398 &   0.9740 &   0.9899\\ 
$2^{-20}$ &   0.7703 &   0.5813 &   0.7975  &   0.8797  &   0.9398 &   0.9740 &   0.9898\\ \hline 
Uniform &   0.9004 &   0.5813 &   0.7975  &   0.8797 &   0.9398 &   0.9740  &   0.9898\\  \hline
\end{tabular}

\end{table}

\begin{table}\label{T6}
\caption{Orders of global convergence for upwinding applied to  problem  (\ref{test1})} 
\begin{tabular}{|r|rrrrrrr|}\hline 
$\varepsilon| N$&$N=8$& 16&32&64&128&256&512\\ \hline 
$2^{0}$ &   1.1033 &   1.0507 &   1.0286  &   1.0148  &   1.0075 &   1.0037 &   1.0019\\ 
$2^{-4}$  &   0.7810 &   0.9087 &   0.9019  &   0.9085  &   0.9134 &   0.9193 &   0.9189\\ 
$2^{-8}$  &   0.8370 &   0.8311 &   0.8458  &   0.8674  &   0.8775 &   0.8956 &   0.9074\\ 
$2^{-12}$  &   0.8355 &   0.8268 &   0.8418  &   0.8713  &   0.8752 &   0.8946 &   0.9064\\ 
$2^{-16}$ &   0.8354 &   0.8263 &   0.8419  &   0.8715  &   0.8751 &   0.8945 &   0.9063\\ 
$2^{-20}$ &   0.8354 &   0.8262 &   0.8419  &   0.8715  &   0.8751 &   0.8945 &   0.9063\\ \hline 
Uniform &   0.8354 &   0.8262 &   0.8419  &   0.8715 &   0.8751 &   0.8945  &   0.9063\\ \hline
\end{tabular}
\end{table}

\begin{table}\label{T10}
\caption{Orders of global convergence for the scheme (\ref{test_L^*_lumped})  applied to problem (\ref{test3})}
\begin{tabular}{|r|rrrrrrr|}\hline 
$\varepsilon| N$&$N=8$& 16&32&64&128&256&512\\ \hline 
$2^{0}$&   1.9297 &   1.9604 &   1.9844  &   1.9930  &   1.9967 &   1.9984 &   1.9992\\ 
$2^{-4}$   &   0.5233 &   0.8485 &   1.1258  &   1.3165  &   1.4590 &   1.5645 &   1.6393\\ 
$2^{-8}$  &   0.3879 &   0.6967 &   1.0264  &   1.2831  &   1.4689 &   1.3243 &   1.6966\\ 
$2^{-12}$ &   0.3807 &   0.6921 &   1.0227  &   1.2796  &   1.3505 &   0.9298 &   0.9505\\ 
$2^{-16}$ &   0.3802 &   0.6917 &   1.0224  &   1.2794  &   1.3496 &   0.9295 &   0.9502\\ 
$2^{-20}$&   0.3801 &   0.6917 &   1.0224  &   1.2794  &   1.3496 &   0.9295 &   0.9501\\ \hline 
Uniform &   0.4410 &   0.6917 &   1.0224  &   1.2794 &   1.3496 &   0.9295  &   0.9501\\  \hline
\end{tabular}
\end{table}

\begin{table}\label{T11}
\caption{Orders of global convergence for standard upwinding  applied to the 
 test problem (\ref{test3})}
\begin{tabular}{|r|rrrrrrr|}\hline 
$\varepsilon| N$&$N=8$& 16&32&64&128&256&512\\ \hline 
$2^{0}$   &   0.8903 &   0.9471 &   0.9710  &   0.9847  &   0.9919 &   0.9959 &   0.9979\\ 
$2^{-4}$   &   0.3737 &   0.5625 &   0.6500  &   0.6954  &   0.7152 &   0.8034 &   0.8202\\ 
$2^{-8}$  &   0.3337 &   0.5258 &   0.5879  &   0.6735  &   0.7304 &   0.8087 &   0.8301\\ 
$2^{-12}$  &   0.3299 &   0.5229 &   0.5821  &   0.6726  &   0.7317 &   0.8097 &   0.8303\\ 
$2^{-16}$ &   0.3298 &   0.5227 &   0.5817  &   0.6726  &   0.7318 &   0.8098 &   0.8303\\ 
$2^{-20}$ &   0.3298 &   0.5227 &   0.5817  &   0.6726  &   0.7318 &   0.8098 &   0.8303\\ \hline 
Uniform &   0.3298 &   0.5227 &   0.5817  &   0.6726 &   0.7318 &   0.8098  &   0.8303\\  \hline
\end{tabular}
\end{table}

\end{document}